\documentclass[12pt]{article}

\usepackage{diagram,amssymb}

\newcommand{\ie}{{\it i.e.\ }}
\newcommand{\eg}{{\it e.g.\ }}

\newcommand{\Cm}{\mathbb{C}}
\renewcommand{\H}{\mathcal{H}}

\newcommand{\N}{\mathcal{N}}
\newcommand{\R}{\mathbb{R}}
\newcommand{\So}{\mathcal{S}}

\newcommand{\bbar}{\bar{\beta}}
\newcommand{\fbar}{\bar{f}}
\newcommand{\thetatil}{\tilde{\theta}}

\newcommand{\Gkn}{\ensuremath{G_{k,n}}}
\newcommand{\Ker}{\mathrm{Ker\,}}
\newcommand{\Skn}{\ensuremath{{\bf S}_{k,n}}}
\newcommand{\Qkn}{\ensuremath{{\bf Q}_{k,n}}}
\newcommand{\Tkn}{\ensuremath{{\bf T}_{k,n}}}
\newcommand{\Z}{\mathbb{Z}}

\newtheorem{thm}[subsection]{Theorem}
\newtheorem{prop}[subsection]{Proposition}
\newenvironment{prf}{{\noindent \hspace{1.3em} \bf Proof.}}
                     {\hspace*{\fill} $\blacksquare$ \vspace{4pt}}

\begin{document}

\title{The quantum Euler class and the quantum cohomology of the
Grassmannians}
\author{Lowell Abrams}

\maketitle

\begin{abstract}
The Poincar\'{e} duality of classical cohomology and the extension of
this duality to quantum cohomology endows these rings with the
structure of a Frobenius algebra. Any such algebra possesses a
canonical ``characteristic element;'' in the classical case this is
the Euler class, and in the quantum case this is a deformation of the
classical Euler class which we call the ``quantum Euler class.'' We
prove that the characteristic element of a Frobenius algebra $A$ is a
unit if and only if $A$ is semisimple, and then apply this result to
the cases of the quantum cohomology of the finite complex
Grassmannians, and to the quantum cohomology of hypersurfaces. In
addition we show that, in the case of the Grassmannians, the [quantum]
Euler class equals, as [quantum] cohomology element and up to sign,
the determinant of the Hessian of the [quantum] Landau-Ginzbug
potential.
\end{abstract}

\section{Introduction}

In \cite{Wit82}, Witten's study of instantons in the context of
supersymmetry of systems with deformed Hamiltonians gave rise to the
notion of a deformed cohomology ring. This ``quantum cohomology ring''
has since then been formulated precisely in terms of Gromov-Witten
invariants of symplectic manifolds (see \cite{McDSal} for
details). Necessarily, much of the attention paid to quantum
cohomology has been from the point of view of symplectic geometry, \eg
\cite{RuaTia95, McDSal}. There has also been a great deal of natural
interest in the realm of algebraic geometry, \eg \cite{KonMan94,
CraMir94}.

Nevertheless, there is strong motivation to pursue an approach which
emphasizes and investigates the parallels between classical and
quantum cohomology.  The quantum cohomology ring of a manifold $M$ is
additively essentially the same as the classical cohomology ring of
$M$, but possesses a multiplication which is a deformation of the
classical cup product (see section \ref{secqcoh}). The strong analogy
between the algebraic structures of these two rings is responsible for
the fact that the Euler class has a quantum analogue which we refer to
as the ``quantum Euler class,'' defined in $\S 1$ below. We show here
that this element of the quantum cohomology ring carries with it
information about the semisimplicity, or lack therof, of the quantum
cohomology ring. 

The issue of semisimplicity of quantum cohomology rings has already
been under investigation from other points of view, as in \cite{Dub96,
KonMan94}. In \cite{Dub96}), Dubrovin defines a {\bf Frobenius
manifold} $M$ to be a manifold such that each fiber of the tangent
bundle $TM$ has a Frobenius algebra (FA) structure, which varies
``nicely'' from fiber to fiber. This context allows for a close
investigation of the nature of the quantum deformations of classical
cohomology, which is generally realized as $T_0M$, the tangent plane
at ``the origin'' in $M$. Moreover, the fact that $M$ is a Frobenius
manifold is equivalent to the existence of a ``Gromov-Witten
potential'' on $M$ satisfying various differential equations,
including the ``WDVV'' equations [{\it ibid}, p. 133].  Special
manifolds, which Dubrovin calls {\bf massive Frobenius manifolds},
have the additional property that for a generic point $t \in M$, the
FA $T_tM$ is semisimple. In this case, a variety of additional results
relating to the classification of Frobenius manifolds hold [{\it
ibid}, Lecture 3].

Kontsevich and Manin discuss aspects of Frobenius manifolds in
\cite{KonMan94}, but deal with a different notion of
semisimplicity. Working with a manifold $M$ which is essentially the
cohomology ring of some space, they define a particular section $K
\colon M \rightarrow TM$ and, at each point $\gamma \in M$, the linear
operator $B(\gamma) \colon T_{\gamma}M \rightarrow T_{\gamma}M$ which
is ``multiplication by $K(\gamma)$.'' They also define a particular
extension $\tilde{T}M$ of $TM$ and show that if, over a subdomain of
$M$, the operator $B(\gamma)$ is semisimple (\ie has distinct
eigenvalues), then $\tilde{T}M$ exhibits some special properties. The
notion of semisimplicity of $B(\gamma)$ is referred to as
``semisimplicity in the sense of Dubrovin'' in \cite{TiaXu96} and
other locations.

The quantum Euler class defined here also provides a section of the
tangent bundle of a Frobenius manifold, although its exact connection
with semisimplicity in the sense of Dubrovin is not yet clear.

The general structure and content of this article are as follows: The
expository presentation of classical cohomology in \S
\ref{secclasscoh} highlights the algebraic structures which are
generalized and deformed in \S\S \ref{secFA} - \ref{secqcohhyp}. In
particular, we offer a new canonical description of the Euler class
$e$. The approach of \S \ref{secclasscoh} is extended to the general
case of Frobenius algebras in \S \ref{secFA}, where the generalized
analogue of the Euler class -- ``the characteristic element'' -- is
shown to satisfy the following:

\vspace{.5\baselineskip}
\noindent 
{\bf Theorem \ref{thomegaunit}} \ {\it The characteristic element of a
Frobenius algebra $A$ is a unit if and only if $A$ is semisimple.}
\vspace{.5\baselineskip}
 
Strictly speaking, quantum cohomology should be viewed as a ring
extension, and not an algebra. Section \ref{secFE} provides the
algebraic framework necessary to generalize the material of \S
\ref{secFA} to the case of a Frobenius extension (FE), \ie when the
base ring is not a field. This having been done, \S \ref{secqcoh}
sketches the elements of the definition of quantum cohomology,
emphasizing its structure as a deformation of classical cohomology,
and in particular as a FE. The ``quantum Euler class'' $e_q$, which is
a deformation of $e$, is defined here to be the characteristic element
of the FE structure of the quantum cohomology ring. Utilizing the
material in \S \ref{secFE}, the semisimplicity test \ref{thomegaunit}
can be applied to quantum cohomology rings.

In the classical and quantum cohomology rings of the complex
Grassmannians, the Euler class and quantum Euler class take on
additional significance. Section \ref{secqcohgrass} outlines how these
rings can be described as Jacobian algebras, where the ideal of
relations is generated by the partial derivatives of the appropriate
Landau-Ginzburg potential $W$ ($W_q$ in the quantum case). In this
context we prove the following result:

\vspace{.5\baselineskip}
\noindent 
{\bf Theorem \ref{thhesseuler}} \ {\it The classical and quantum Euler
classes are equal, up to sign, to the determinants of the Hessians of
$W$ and $W_q$, respectively.}
\vspace{.5\baselineskip}

\noindent
This connects the classical and quantum Euler classes to
Morse-theoretic considerations regarding the functions $W$ and
$W_q$. In a sense, it brings them back to the roots of quantum
cohomology in \cite{Wit82}, which utilizes a Morse-theoretic
approach. In addition, this result leads to a new proof of proposition
\ref{prqcohsemi} which, modulo technicalities, states that the quantum
cohomology of any finite complex Grassmannian manifold is semisimple.

Finally, \S \ref{secqcohhyp} applies the semisimplicity test
\ref{thomegaunit} to the quantum cohomology of hyperplanes, providing
good contrast to the situation for the Grassmannians.

\section{Classical Cohomology and the Euler Class}
   \label{secclasscoh}

Let $X$ denote a connected $K$-oriented $n$-dimensional compact
manifold, where $n$ is even. Throughout this article, except where
noted otherwise, homology and cohomology groups will use coefficients
in a field $K$ of characteristic $0$. Denote by $[X] \in H_n(X)$ the
fundamental orientation class of $X$, and let $\langle-,-\rangle
\colon H^*(X) \otimes H_*(X) \rightarrow K$ denote the Kronecker
index. The kernel of the linear form $\mu^* \colon H^*(X) \rightarrow
K$ , where $\mu$ denotes the generator of $H^n(X)$ satisfying $\langle
\mu, [X] \rangle = 1$, contains no nontrivial ideals. This form can be
used to define the ``intersection form'' $H^*(X) \otimes H^*(X)
\rightarrow K$, by $a \otimes b \mapsto \mu*(a \cup b)$. The
intersection form is nondegenerate.

Notice that we may view $H_*(X)$ as a (left) $H^*(X)$-module via the
cap product $\cap \colon H^*(X) \otimes H_*(X) \rightarrow
H_*(X)$. Viewing $H^*(X)$ as the regular (left) module over itself, we
see that the Poincar\'{e} duality map
\[
D \colon H^*(X) \rightarrow H_*(X), \ \ \zeta \mapsto \langle-,\zeta
   \cap [X]\rangle
\]
is an $H^*(X)$-module isomorphism.

Let $\Delta \colon X \rightarrow X \times X$ denote the diagonal
map. The transfer map $\Delta^! \colon H^*(X) \rightarrow H^*(X)
\otimes H^*(X)$ is defined to be the map which makes the following
diagram commutative:
\[
\begin{diagram}
   \node{H^*(X)} \arrow{e,t}{\Delta^!} \arrow{s,l}{D} \node{H^*(X)
      \otimes H^*(X)} \\ \node{H_*(X)} \arrow{e,t}{\Delta_*}
      \node{H_*(X) \otimes H_*(X)} \arrow{n,l}{D^{-1} \otimes D^{-1}}
\end{diagram}
\]
Here, we implicitly use the isomorphism $H_*(X \times X) \cong H_*(X)
\otimes H_*(X)$, and the corresponding isomorphism for
cohomology. Modulo this latter isomorphism, the cup-product in
$H^*(X)$ is given by $\Delta^* \colon H^*(X) \otimes H^*(X)
\rightarrow H^*(X)$.

Let $j \colon (X \times X, \emptyset) \rightarrow (X \times X,\, X
\times X \setminus \Delta(X))$ denote inclusion of pairs. Consider the
element $\tau := \Delta_!(1) = (D^{-1} \otimes D^{-1}) \circ
\Delta_*([X])$. By the canonical isomorphism of the tangent bundle
$TX$ to the normal bundle of $\Delta(X)$ in $X \times X$
\cite{MilSta}, this is just the image under $j^*$ of the Thom class of
$TX$. It follows that $\Delta^* \circ \Delta^!(1) \in H^*(X)$ is in
fact the Euler class $e(X)$.

We recall the well known formula \cite{MilSta}
\[
e(X) \ = \ \sum_i e_ie_i^{\#},
\]
where $e_i$ ranges over a basis for $H^*(X)$, and $e^{\#}_j$ ranges
over the corresponding dual basis relative to the intersection form,
\ie $\mu^*(e_i \cup e_j^{\#}) = \delta_{ij}$.

\section{Frobenius Algebras and the Characteristic Element}
   \label{secFA}

Let $K$ be a field of characteristic $0$ and let $A$ be a
finite-dimensional (as a vector space) commutative algebra over $K$,
with unity $1_A$. Let $\beta \colon A \otimes A \rightarrow A$ denote
multiplication in $A$, and let $\bbar \colon A \rightarrow
\mathrm{End}(A)$ denote the regular representation of $A$, \ie
$\bbar(a)$ is ``multiplication by $a$.'' View $A$ as the regular
module over itself, and view the vector space dual $A^*$ as an
$A$-module via the action $A \otimes A^* \rightarrow A^*$ given by $a
\otimes \zeta \, \mapsto \, a \cdot \zeta := \zeta \circ \bbar(a)$.

$A$ is referred to as a {\bf Frobenius algebra} (FA) if there exists
an $A$-module isomorphism $\lambda \colon A \rightarrow A^*$, \ie a
nondegenerate pairing. In \cite[pages 414-418]{CurRei} this is shown
to be equivalent to the existence of a linear form $f \colon A
\rightarrow K$ whose kernel contains no nontrivial ideals, and to the
existence of a nondegenerate linear form $\eta \colon A \otimes A
\rightarrow K$ which is associative, \ie $\eta(ab \otimes c) = \eta(a
\otimes bc)$. In fact, we may take $f := \lambda(1_A)$ and $\eta := f
\circ \beta$, and we will henceforth presume that $\lambda, f$ and
$\eta$ are related in this way. When it is useful to emphasize the FA
structure of $A$ endowed by particular $f, \eta,$ and $\lambda$, the
algebra $A$ will be denoted by $(A,f)$.

For the next result, view $A \otimes A$ as an $A$-module via the usual
module action $\beta \otimes I \colon A \otimes A \otimes A
\rightarrow A \otimes A$.

\begin{thm}
   \label{thFAcoalg} A finite dimensional commutative algebra $A$ with
$1_A$ is a FA if and only if it has a cocommutative comultiplication
$\alpha \colon A \rightarrow A \otimes A$, with a counit, which is a
map of $A$-modules.
\end{thm}

\begin{prf}
A complete proof appears in \cite{Abr96}. Here, we simply note that if
$A$ is a Frobenius algebra with pairing $\lambda$, then the
comultiplication $\alpha$ is defined to be the map $(\lambda^{-1}
\otimes \lambda^{-1}) \circ \beta^* \circ \lambda$:
\[
\begin{diagram}
   \node{A} \arrow{e,t}{\alpha} \arrow{s,l}{\lambda} \node{A \otimes
      A}\\ \node{A^*} \arrow{e,t}{\beta ^*} \node{A^* \otimes A^*}
      \arrow{n,r}{\lambda^{-1} \otimes \lambda^{-1}}
\end{diagram}
\]
\enlargethispage*{10\baselineskip}
\end{prf}
\pagebreak

Define the {\bf characteristic element of $(A,f)$} to be the element $
\omega_{A,f} := \beta \circ \alpha(1_A) \in A. $ This is a canonical
element which is shown in \cite{Abr96} to be of the form
\[
\omega_{A,f} \, = \, \sum_i e_ie^{\#}_i \, ,
\]
where $e_i$ ranges over a basis for $A$ and $e_j^{\#}$ ranges over the
corresponding dual basis relative to $\eta$.

It is easy to show that theorem \ref{thFAcoalg} still holds if
``commutative'' is replaced by ``skew-commutative,'' as would be the
case for $H^*(X)$. We see that in that case $f, \lambda, \alpha,
\omega$ correspond to $\mu^*, D, \Delta^!, e(X)$, respectively.

Given FA's $(A,f)$ and $(B,g)$, we can form the {\bf direct sum} $(A
\oplus B, f \oplus g)$, where $A \oplus B$ denotes the ``orthogonal
direct sum'' of algebras, and $f \oplus g$ acts by $f \oplus g(a
\oplus b) := f(a)+g(b) \in K$.  $(A \oplus B, f \oplus g)$ is in fact
a FA \cite{Abr96}.

\begin{prop} \mbox{\rm \cite{Abr96}}
   \label{promegarespect} The characteristic element respects direct
 sum structure. Specifically,
\[
\omega_{A' \oplus A'' , f' \oplus f''} = \omega_{A',f'} \oplus
\omega_{A'',f''} \in A' \oplus A''
\]
\end{prop}

The minimal essential ideal $\So = \So(A)$ of a ring $A$ is called the
{\bf socle}. When $A$ is indecomposable, the socle is
$\mathrm{ann}({\N})$, where $\N = \N(A) \subset A$ is the ideal of
nilpotents. See \cite[\S 9]{AndFul} for details.

\begin{prop}
   \label{promegasocle} In a FA $A$, the ideal $\omega A$ is the socle
of A.
\end{prop}

This result is independent of the choice of FA structure.

The construction in the following proof is essentially taken from
Sawin \cite{Saw95}, although this result does not explicitly appear
there.

\begin{prf}
Because the socle of a finite-dimensional commutative algebra is the
direct sum of the socles of its indecomposable constituents \cite[\S
9]{AndFul}, it suffices to prove this proposition for the
indecomposable cases. Furthermore, we showed in \cite{Abr96} that the
socle $\So$ of a FA is a principal ideal, any of whose elements is a
generator, so it suffices to show that $\omega$ lies in the
socle. Notice that $\omega$ is not $0$; we have $f(\omega) = (A : K)
\in K$, and this is not $0$ in $K$, since $K$ has characteristic $0$.

If $A$ is a field extension then $\N(A) = \{0\}$, so the socle $\So =
\mathrm{ann}(\N) = A$. But $\omega$ is not zero, so it is a unit, and
thus $\omega A = A = \So$.

If $A$ is not a field extension, define a chain of ideals $\So=S_1
\subset S_2 \subset \cdots \subset S_n=A$, where each $S_k$ is the
preimage in $A$ of the socle of $A/S_{k-1}$. Choose a basis for
$S_1$. Now, starting with $i=1$, iteratively take the basis for $S_i$
and extend it to a basis for $S_{i+1}$. Denote the elements of the
basis for $S_n=A$ by $e_1, \ldots , e_n$, and let $e_1^{\#}, \ldots ,
e_n^{\#}$ denote the corresponding dual basis elements. Suppose $e_i
\in S_k \setminus S_{k-1}$ and that $a \in A$ is any nilpotent
element. Then $ae_i \in S_{k-1}$, and therefore can be expressed as a
linear combination of basis elements other than $e_i$. It follows that
$f(ae_ie_i^{\#}) = 0$, so $e_ie_i^{\#}\N(A) \subset \Ker f$. But $\Ker
f$ can contain no nontrivial ideals, as mentioned above, so we must
have $e_ie_i^{\#}\N(A) = 0$, ie. $e_ie_i^{\#} \in \So$. This follows
for each $i$, so $\omega = \sum_i e_ie_i^{\#} \in \So$.
\end{prf}

\begin{thm}
   \label{thomegaunit} The characteristic element $\omega$ of a FA $A$
is a unit if and only if $A$ is semisimple.
\end{thm}
\begin{prf}
First, recall from the proof of \ref{promegasocle} that $\omega$ is
not $0$. Because $A$ is commutative, it is semisimple if and only if
it is a direct sum of fields. In such a case, the component of
$\omega$ in each component of $A$ is nonzero (each component is a FA
\cite{Abr96}), and hence a unit. Since a direct sum of units is a
unit, $\omega$ is a unit.

If some component $A'$ of $A$ is not a field, then it contains
nontrivial nilpotents. In this case, $\So(A') = \mathrm{ann}(\N(A'))$
is nilpotent, so $\omega$ has a nilpotent component, and cannot be a
unit.
\end{prf}

In a skew-commutative context, such as $H^*(X)$, the characteristic
element is not necessarily nonzero. For instance, if $X$ is an
odd-dimensional compact oriented manifold then the characteristic
element, \ie the Euler class, is $0$. However, if the characteristic
element is in fact nonzero, then \ref{thomegaunit} still holds.

\section{Frobenius Extensions}
   \label{secFE}

Suppose $A/R$ is a finite-dimensional (as a module) commutative ring
extension with identity. By analogy with FA's, if there exists a
module isomorphism $\lambda \colon A \rightarrow A^*$, we call $A$ a
{\bf Frobenius extension} (FE). As in the case of FA's, this is
equivalent to the existence of maps $\eta$ and $\alpha$. There is also
a ``FE form'' $f := \lambda(1_A) \colon A \rightarrow R$, but in this
context it is not sufficient for the kernel of $f$ to contain no
nontrivial ideals. The characteristic element $\omega_{A,f}$ may be
defined as for FA's, but note that theorem \ref{thomegaunit} no longer
applies. This section provides an approach for dealing with this
circumstance.

Suppose $\theta \colon R \rightarrow S$ is a surjective homomorphism
of rings (sending $1_R \mapsto 1_S$). Let $(A,f)$ denote a FE, and
define $B = \theta_*(A)$ to be $A \otimes_R S$. In this ring, we have
$ra \otimes s \, = \, a \otimes \theta(r)s$ for all $r \in R,\, s \in
S$, and $a \in A$. Let $\thetatil \colon A \rightarrow B$ denote the
ring homomorphism $a \mapsto a \otimes 1_S$. Define the linear form
$\fbar \colon B \rightarrow S$ by
\[
\fbar(a \otimes s) := \theta \circ f(a)s.
\]
The form $\fbar$ is well-defined, since
\[
\fbar(ra \otimes s) \, = \, \theta \circ f(ra)s \, = \, \theta(rf(a))s
                    \, = \, \theta(r)(\theta \circ f(a))s \, = \,
                    \fbar(a \otimes \theta(r)s),
\]
and $\fbar$ satisfies the commutative diagram \[
\begin{diagram}
   \node{A} \arrow{e,t}{\thetatil} \arrow{s,r}{f} \node{B}
      \arrow{s,r}{\fbar} \\ \node{R} \arrow{e,t}{\theta} \node{S}
\end{diagram}
\]

Let $e_1, \ldots, e_n$ denote a basis for $A$, and let $e_i^{\#},
\ldots, e_n^{\#}$ denote the corresponding dual basis relative to
$\eta_A$.

\begin{prop}
   \label{prinducedFE} The form $\fbar$ endows $B=\theta_*(A)$ with a
FE structure, and
\[
\omega_{B,\fbar} = \thetatil(\omega_{A,f}).
\]
\end{prop}

\begin{prf}
It suffices to show that the set $\{\thetatil(e_1), \ldots,
\thetatil(e_n)\}$ is a basis for $B$, and that its dual basis relative
to the form $B \otimes B \rightarrow S$, $a \otimes b \mapsto
\fbar(ab)$ is $\{\thetatil(e_1^{\#}), \ldots,
\thetatil(e_n^{\#})\}$. The existence of a dual basis will show $B$ is
a FE. The particular form of the basis and dual basis, together with
the fact that $\thetatil$ is a homomorphism, will prove the claim
about $\omega_{B,\fbar}$.

We first prove the orthogonality relations:
\[
\fbar \left( \thetatil(e_i)\thetatil(e_j^{\#}) \right) = \fbar \left(
  \thetatil(e_ie_i^{\#}) \right) = \theta \circ f(e_ie_i^{\#}) =
  \theta(\delta_{ij}) = \delta_{ij} \in S.
\]
To prove that we have a basis as claimed, note that the elements
$\thetatil(e_1), \ldots, \thetatil(e_n)$ clearly span $B$, since
$\thetatil$ is surjective. Suppose that for some $\{s_i\} \in S$ we
have $\sum_i s_i\thetatil(e_i) = 0$. Then, for all $j$,
\[
0 = \fbar(0) = \fbar \left( \sum_i
  s_i\thetatil(e_i)\thetatil(e_j^{\#}) \right) = s_j.
\]
It follows that $\thetatil(e_1), \ldots, \thetatil(e_n)$ are
independent, and thus form a basis. The orthogonality relations show
that $\{\thetatil(e_1^{\#}), \ldots, \thetatil(e_n^{\#})\}$ is a basis
as well.
\end{prf}

In the next result, let $\theta \colon R \rightarrow K$ be any
surjective $K$-linear ring homomorphism, where $K$ is a field.

\begin{prop} \hspace{2pt}
   \label{prinducedomega} The element $\omega_{A,f}$ is either a unit
in $A$ or a zero divisor.
\begin{description}
   \item[(i)] If $\omega_{A,f}$ is a unit in $A$ then $B =
   \theta_*(A)$ is semisimple.  \item[(ii)] If $\omega_{A,f}$ is a
   zero divisor and \ $\mathrm{ann}(\omega_{A,f}) \nsubseteq \Ker
   \thetatil$, then $\theta_*(A)$ is not semisimple.
\end{description}
\end{prop}

\begin{prf}
If $\omega_{A,f}$ is a unit, then there exists a $u \in A$ such that
$\omega_{A,f} u = 1_A$. But then, by \ref{prinducedFE}
\[
\omega_{B,\fbar}\thetatil(u) = \thetatil(\omega_{A,f})\thetatil(u) =
\thetatil(1_A) = 1_B,
\]
so $\omega_{B,f}$ is a unit as well.

All FE structures on $A$ are given by $(A,f \circ \bbar(u))$, for some
unit $u \in A$ \cite[Proposition 2, {\it mutatis
mutandis}]{Abr96}. Thus, if $\omega$ is not a unit in $A$, then the
map $\omega_{A,f} \cdot f$ is not a FE form. This implies that there
exists an $a \in A$ such that $f(\omega_{A,f} aA) = \omega_{A,f} \cdot
f(aA) = \{0\}$. But $f$ is a FE form, so it must be that $\omega_{A,f}
a = 0$. If follows that $\thetatil(\omega_{A,f})\thetatil(a) =
0$. Since, by assumption, there exists some $a \in
\mathrm{ann}(\omega_{A,f})$ such that $a \notin \Ker \thetatil$, we
see that $\thetatil(\omega_{A,f}) = \omega_{B,\fbar}$ is a zero
divisor as well. Both statements (i) and (ii) now follow from theorem
\ref{thomegaunit}.
\end{prf}

\section{Quantum Cohomology and the Quantum Euler Class}
   \label{secqcoh}

Let $X$ be a $2n$-dimensional compact oriented manifold which, in
addition, is symplectic, and let $H'_2(X)$ denote the free part of
$H_2(X,\Z)$. Taking $B_1, \ldots, B_n$ to denote a basis of $H'_2(X)$,
the group algebra $\Lambda := K[H'_2(X)]$ may be expressed as
$K[q^{B_1}, \ldots, q^{B_n}]$, where $q$ is a formal variable and the
addition of exponents is the group operation of $H'_2(X)$. This is
essentially an algebraic version of the Novikov ring (see \cite[\S
9.2]{McDSal}). As an additive group, the {\bf quantum cohomology} ring
$\mathbf{QH^*(X)}$ has the same structure as $H^*(X) \otimes \Lambda$,
but has a ``deformed'' multiplication, which we describe briefly:

The classical cup product of two elements $a, b \in H^*(X)$ is given
by
\[
a \cup b = \sum_i(\alpha \cdot \beta \cdot \gamma_i)c_i \, ,
\]
where $c_i$ runs over a basis for $H^*(X)$ and $\alpha, \beta,
\gamma_i$ are the Poincar\'{e} duals of $a, b, c_i$, respectively, and
``$ \cdot$'' denotes the homology intersection index. The quantum
multiplication
\[
\ast \colon QH^*(X) \otimes QH^*(X) \rightarrow QH^*(X)
\]
is defined on elements $a,b \in H^*(X) \hookrightarrow QH^*(X)$ by
\[
a \ast b \, := \, \sum_{i,B} \Phi_B(\alpha, \beta, \gamma_i)q^Bc_i \,
,
\]
and extended by linearity to all of $QH^*(X)$.  Here, $B$ ranges over
$H'_2(X)$, and $\Phi_B(\alpha, \beta, \gamma_i)$ denotes the Gromov
(Gromov-Witten) invariants. Intuitively, these count intersections
(subject to dimension requirements!) of the cells $\alpha, \beta,
\gamma_i$ not with themselves, but with the fourth cell $B$. When $B =
0$, the Gromov invariant is the classical intersection index. Thus,
\[
a \ast b \ = \ a \cup b \, + \mathrm{other\ terms}.
\]
For details regarding the definition of quantum cohomology, and in
particular proofs of the associativity of $\ast$, see
\cite{McDSal,RuaTia95}.

Extend $\mu^* \colon H^*(X) \rightarrow K$ (defined in section
\ref{secclasscoh}) by linearity over $\Lambda$ to a form $\mu^* \colon
QH^*(X) \rightarrow \Lambda$.
\begin{prop}
   \label{prqcohFE} The form $\mu^*$ endows $QH^*(X)$ with a FE
structure.
\end{prop}
\begin{prf}
   See \cite{Abrth} for a rigorous proof.
\end{prf}

Let $\iota \colon H^*(X) \hookrightarrow QH^*(X)$ denote the obvious
inclusion map.  Note that although $H^*(X)$ and $QH^*(X)$ share
essentially the same basis $\{e_i\}$ and the same FA form, the
respective dual bases are not necessarily equal. In other words, the
fact that the element $e_i^{\#}$ is the dual in $H^*(X)$ to $e_i$ does
not necessarily imply that $\iota(e_i^{\#})$ is dual to $\iota(e_i)$
in $QH^*(X)$. However, it does hold that the $q^0$ term of
$\iota(e_i)^{\#}$ is in fact $\iota(e_i^{\#})$. It follows that the
$q^0$ term of the characteristic element $\omega_q$ of $(QH^*(X),
\mu^*)$ is $e(X)$. In other words, we have:

\begin{center} 
   The characteristic element $\omega_q$ is a deformation of the
classical Euler class.
\end{center}

Because of this, we refer to $\omega_q$ as {\bf the quantum Euler
class}, and denote it by $e_q(X)$. Unlike $e(X)$, the quantum Euler
class may very well be a unit. Strictly speaking, however, the
semisimplicity result \ref{thomegaunit} does not apply to $QH^*(X)$,
because it is infinite dimensional (as a vector space) over $K$ and
only a ring extension (not an algebra) over $\Lambda.$ We may,
however, utilize the approach of section \ref{secFE}. Define the
homomorphism $\theta \colon \Lambda \rightarrow K$ as follows: For
each generator $B_i$ of $H'_2(X)$ choose any nonzero $r_i \in K$ and
define $\theta(q^{B_i}) := r_i$. Extending $\theta$ by linearity over
$K$ gives a surjective ring homomorphism, often referred to as
``specialization.'' Theorem \ref{thomegaunit} now applies to
$\theta_*[QH^*(X)]$, which is a FA.

\section{The Quantum Cohomology of the Grassmannians}
   \label{secqcohgrass}

Let \Gkn\ denote the Grassmannian manifold of complex $k$-dimensional
subspaces in $\Cm^n$. Define the Chern polynomial of $X = \Gkn$ to be
\[
c_t(\Gkn) := \sum_{i=1}^k x_it^i = \prod_{i=1}^k(1 + \lambda_it),
\]
where $t$ is a formal variable and the $x_i$'s are the Chern classes
of the canonical bundle \Skn. The $\lambda_i$ are referred to as the
{\bf Chern roots} of \Gkn\ (but they are {\it not} roots of
$c_t$!). Obviously, $x_i$ is the $i$'th elementary symmetric
polynomial $\sigma_i(\lambda_1, \ldots, \lambda_k)$ in the Chern
roots. Define
\[
W(\lambda_1, \ldots, \lambda_k) = \sum_{i=1}^k
\frac{1}{n+1}\lambda_i^{n+1}
\]
\begin{center} and \end{center}
\[
\begin{array}{rcl}
W_q(\lambda_1, \ldots, \lambda_k) & = & \displaystyle{ \sum_{i=1}^k
   \frac{1}{n+1}\lambda_i^{n+1} + (-1)^kq\lambda_i} \\ \\ & = &
   W(\lambda_1, \ldots, \lambda_k) + (-1)^kqx_1.
\end{array}
\]
The function $W_q$ is called the Landau-Ginzburg potential of
\Gkn. Because $W$ and $W_q$ are symmetric functions in the
$\lambda_i$, they may also be viewed as functions of $x_1, \ldots,
x_k$. Define $dW$ to be the ideal $( \frac{\partial W}{\partial x_1},
\ldots, \frac{\partial W}{\partial x_k} )$, and define $dW_q$
similarly. Then
\[
H^*(\Gkn) \cong K[x_1, \ldots, x_k] / dW
\]
\begin{center} and \end{center}
\[
QH^*(\Gkn) \cong K[q,q^{-1}][x_1, \ldots, x_k] / dW_q .
\]
Denote by $H$ and $H_q$ the determinants of the Hessians
\[
\H = \left( \frac{\partial^2 W}{\partial x_i \partial x_j} \right) \ ,
\ \ \H_q = \left( \frac{\partial^2 W_q}{\partial x_i \partial x_j}
\right)
\]
of $W$ and $W_q$, respectively.

In this section, we will prove the following:
\begin{thm}
   \label{thhesseuler} $e(\Gkn) = (-1)^{n \choose 2}H$ \ and \
$e_q(\Gkn) = (-1)^{n \choose 2}H_q$.
\end{thm}

Suppose that an algebra $A$ (not necessarily a FA) is finite
dimensional as a vector space and is given by the presentation $A
\cong K[x_1, \ldots, x_n]/R$, where $R = (f_1, \ldots, f_p)$ is some
finitely-generated ideal in $K[x_1, \ldots, x_n]$. Note that we
continue to assume that $K$ has characteristic $0$. Because $A$ is
finite dimensional, we must have $p \geq n$. The {\bf Jacobian ideal}
$J = J(R)$ of $R$ is defined to be the ideal generated by the
determinants of the $n \times n$ minors of the matrix
\[
\left( \frac{\partial(f_1, \ldots, f_p)}{\partial(x_1, \ldots, x_n)}
\right) \mathrm{mod}\, R.
\]
The ideal $J$ is well-defined since it is a Fitting ideal of the
module $\Omega_{A/K}$ of K\"{a}hler differentials of $A$ (see \cite[\S
1.1, \S 10.3]{Vas}).

The following result of Scheja and Storch \cite{SchSto75} is reported
in more generality in \cite[{\it ibid}]{Vas} (although for the
definition of ``complete intersection'' we refer the reader to
\cite{Kun}):

\begin{prop}
   \label{prjacobiansocle} $J \neq \{0\}$ if and only if $A$ is a
complete intersection and $J$ generates the socle of $A$.
\end{prop}

Now assume that $A \cong K[x_1, \ldots, x_n]/R$ is a FA with
characteristic element $\omega$ for some choice of FA structure.
\begin{prop}
   \label{promegajacobian} $J \neq \{0\}$ if and only if $J = \omega
A$. If $p=n$, then $J \neq \{0\}$ if and only if
\[
\mathrm{det} \left( \frac{\partial f_i}{\partial x_j} \right)
\mathrm{mod}\, R = u \omega
\]
for some unit $u \in A$.
\end{prop}

\begin{prf}
This proposition follows immediately from \ref{promegasocle} and
\ref{prjacobiansocle}.
\end{prf}

\begin{prop}
   \label{prHessomega} For each \Gkn\ there is a $\kappa \in K$ such
that $H = \kappa e(\Gkn)$ and $H_q = \kappa e_q(\Gkn)$.
\end{prop}
\begin{prf}
Because $H$ and $e(\Gkn)$ are the $q^0$ terms of $H_q$ and
$e_q(\Gkn)$, respectively, it suffices to prove the proposition for
the quantum case.

The polynomial $W_q$ is homogeneous of degree $2(n+1)$ \cite[\S
8.4]{McDSal}. In other words, each summand of $W_q$ has degree
$2(n+1)$ in $QH^*(\Gkn)$, where $q$ is taken to have degree
$2n$. Also, $|x_i| = 2i$ for each $i$. Thus, for fixed $i, j$ we have
\[
\left| \left( \frac{\partial^2 W_q}{\partial x_i \partial x_j} \right)
   \right| \ = \ |W_q| - |x_i| - |x_j| \ = \ |W_q| - 2i - 2j.
\]

We now show by induction that $H_q$ is homogeneous of degree
$2k(n-k)$. Each $s \times s$ minor $M$ of $\H$ is a matrix with
entries $m_{ij} := \H_{ij}$ where $i$ and $j$ run over elements of
some ordered subsets $I, J \subset \{1, \dots, k\}$ respectively, and
$\#I = \#J = s$. Define $M\widehat{(i,j)}$ to be the minor of $M$
which does not include the entry $m_{ij}$. We have already shown that
when $s=1$, the single entry of each $M$ is homogeneous of degree
$|W_q| - 2i - 2j$. Assume that for all minors $M$ of $\H$ of size less
than $(s+1) \times (s+1)$, the determinant of $M$ is homogeneous of
degree
\[
s |W_q| - \, 2\sum_{i \in I}i \ - \ 2\sum_{j \in J}j.
\]
Now consider any $(s+1) \times (s+1)$ minor $M$ of $\H$ with index
sets $I, J$. Take any $i' \in I$. Then
\[
\mathrm{det} M = \sum_{j' \in J}
    \mathrm{sgn}(i',j')m_{i'j'}\mathrm{det}M\widehat{(i',j')},
\]
where $\mathrm{sgn}$ is the appropriate function $I \times J
\rightarrow \{+1, -1\}$. By the induction hypothesis,
\[
\begin{array}{rcl}
   \left| m_{i'j'}\mathrm{det}M\widehat{(i',j')} \right| & = &
       |m_{i'j'}| + |M\widehat{(i',j')}| \\ \\ & = & |W_q| - 2i' - 2j'
       + s |W_q| - {\displaystyle \left( 2\sum_{i \in
       I\setminus\{i'\}}i \right) - \left( 2\sum_{j \in
       J\setminus\{j'\}}j \right) }\\ \\ & = & (s+1) |W_q| - \
       {\displaystyle 2 \sum_{i \in I}i \ - \ 2\sum_{j \in J}j }.
\end{array}
\]
But this is independent of the choice of $j'$, so det$M$ is
homogeneous of degree $(s+1) |W_q| - \, 2\sum_{i \in I}i \, - \,
2\sum_{j \in J}j$. In particular, we can take $M = \H$, and thus $H$
is homogeneous of degree
\[
k |W_q| \ - \ 2\sum_{i=1}^k i \ - \ 2s\sum_{j=1}^kj \ = \ k(2n+2) -
   2k(k+1) \ = \ 2k(n-k).
\]

Of course, $e_q$ is also homogeneous of degree $2k(n-k)$ since $e_q =
\sum_i e_ie_i^{\#}$, where $e_i$ runs over a basis for $H^*(\Gkn)$,
and since $|e_i^{\#}| = 2k(n-k) - |e_i|$.

Consider the algebra $A:=K[q][x_1,\ldots,x_k]/dW_q$. By the
definitions of $H_q$ and $e_q$, and the nature of the relations given
by $dW_q$, both $H_q$ and $e_q$ may be viewed as elements of $A$. Now,
the proof of \ref{prinducedFE} applies equally well to the algebra
$A$, so it is a FE, and $e_q$ is in fact the characteristic element of
$A$. Proposition \ref{promegajacobian}, shows that $H_q = ve_q(\Gkn)$
for some unit $v \in A$. Of course, $v$ may also be viewed as an
element of $QH^*(\Gkn)$ which simply has no $q^i$-terms with $i<0$.

Write $v = v' + v''$ where $v'$ is homogeneous of degree $0$, and
$v''$ contains no terms of degree $0$. Since $H_q = v'e_q + v''e_q$
and both $H_q$ and $v'e_q$ are homogeneous of degree $2k(n-k)$, we see
that $v''e_q$ must also be homogeneous of this degree. By degree
considerations, we must have $v''e_q = 0$, and thus $H_q = v'e_q$.

Write $v' = \sum_{j \geq 0} v_j q^j$. Since $v'$ is homogeneous of
degree $0$ and $|q^j| = 2jn$, we see that $|v_j| = -2jn$. But $v'$ is
an element of $A$, so we must have $v_j = 0$ for $j \neq 0$. Thus $v'$
may in fact be viewed as a degree $0$ element in $H^*(\Gkn)$. In other
words, $v'$ is an element $\kappa \in K$.
\end{prf}

Take $K=\R$ or $\Cm$, and for any nonzero $r \in K$ let $\theta_r$
denote a specialization homomorphism $K[q,q^{-1}] \rightarrow K,\ q
\mapsto r$ as above. In the following paragraph, any reference to
$QH^*(\Gkn)$ or any element $a$ therein should be interpreted as
referring to $(\theta_r)_*[QH^*(\Gkn)]$ and $\thetatil_r(a)$,
respectively.

In this context, the relationship between the distinguished element
and the Hessian provides $e_q(\Gkn)$ with a nontrivial geometric
interpretation: Denote the critical points of $W_q$ by $z_1, \ldots,
z_j$, and note that $H_q$ may be viewed as a function $K^k \rightarrow
K$, as may all the elements of $QH^*(\Gkn)$. It is well known that,
for each $j$, $H_q(z_j)=0$ if and only if the critical point $z_j$ is
degenerate \cite{Mil}. Because the elements of $QH^*(\Gkn)$, viewed as
functions, are completely determined by their values on the critical
points of $W_q$, we see that $H$ (and hence $e_q(\Gkn)$) is a unit in
$QH^*(\Gkn)$ if and only if the critical points of $W_q$ are all
nondegenerate.

This relationship between $e_q(\Gkn)$ and $H$ also yields a new
approach to the following known result \cite{SieTia94}:

\begin{prop}
   \label{prqcohsemi} For all \Gkn\ and all nonzero $r \in \R$, the
algebra
\begin{center}
$(\theta_r)_*[QH^*(\Gkn)]$ is semisimple.
\end{center}
\end{prop}
The proof is based on calculations appearing in \cite{Ber94}.

\begin{prf}
The Jacobian matrix $V = (\partial x_i / \partial \lambda_j)$
associated to the elementary symmetric functions $x_i$ is a
Vandermonde matrix, and has determinant $\prod_{i<j}(\lambda_i -
\lambda_j) \neq 0$. Let $\nabla_x$ denote the gradient vector operator
with respect to $x_1, \ldots, x_k$, and let $\nabla_{\lambda}$ denote
the gradient operator with respect to $\lambda_1, \ldots,
\lambda_k$. Viewing the gradient operators as row vectors, we have
$\nabla_x(W_q)V = \nabla_{\lambda}(W_q)$. Let $\nabla_x(W_q)_i$ denote
the $i$'th entry of $\nabla_x(W_q)$, and let $V_i$ denote the $i$'th
row of $V$. Then the Hessian of $W_q$ with respect to the $\lambda$'s
is
\[
\begin{array}{rcl}
\nabla_{\lambda}^T \nabla_{\lambda}(W_q) & = & \nabla_{\lambda}^T
     (\nabla_x(W_q)V) \\ \\ & = & V^T \nabla_x^T \nabla_x(W_q) V +
     \sum_i \nabla_x(W_q)_i \, \nabla_{\lambda}^T V_i.
\end{array}
\]
Evaluating at the critical points of $W_q$ (\ie assuming
$\nabla_x(W_q)=0$), and expressing everything in terms of the
$\lambda_i$'s, we see that
\[
H = \mathrm{det} \left( \nabla_{\lambda}^T \nabla_{\lambda}(W_q)
       \right) \mathrm{det} (V^{-2}) = \frac{n^k \prod_{i=1}^k
       \lambda_i^{n-1}} {(\prod_{i<j}(\lambda_i - \lambda_j))^2} .
\]  
Now, because $V$ is invertible, the relation $\nabla_x(W_q)=0$ is
equivalent to $\nabla_{\lambda}(W_q)=0$. In other words, for each $i$
we have $\lambda_i^n = (-1)^{k+1}q$ at the critical points of
$W_q$. This implies that
\[
x_k \prod_{i=1}^k \lambda_i^{n-1} = \prod_{i=1}^k \lambda_i^n =
              (-1)^{k(k+1)}q^k.
\]
Since $q \neq 0$, the numerator of $H$, and thus $H$ itself, is
nonzero at the critical point of $W_q$. It follows that $H$, as an
element of $QH^*(\Gkn)$, has an inverse, and therefore, by
\ref{prHessomega}, so does $e_q(\Gkn)$. By proposition
\ref{prinducedomega}, $\theta_r[QH^*(\Gkn)]$ is semisimple.
\end{prf}

As discussed above, the Chern classes $x_1, \ldots, x_k$ arising from
the bundle \Skn\ are the elementary symmetric polynomials in the Chern
roots $\lambda_1, \ldots, \lambda_k$. An analogous situation holds for
the ``normal'' classes $y_1, \ldots, y_{n-k}$, which arise from the
quotient bundle \Qkn. Define $\mu_1, \ldots, \mu_{n-k}$ to be the
Chern roots corresponding to the formal polynomial
\[
\sum_{i=1}^{n-k} y_it^i.
\]
Then for all $i$, we have $y_i = \sigma_i(\mu_1, \ldots,
\mu_{n-k})$. In fact, the $\lambda_i$'s and $\mu_i$'s are the first
Chern classes of the line bundles in the splitting of \Skn\ and \Qkn,
respectively \cite[\S 21]{BottTu}. Together with the well known
bundle-isomorphism of the tangent bundle $\Tkn \cong \Skn^* \otimes
\Qkn$, this fact allows us to write the characteristic classes
$c_i(\Tkn)$ in terms of the $x_i$'s and $y_i$'s: The Chern polynomial
for \Tkn\ is
\[
\sum_{i=1}^{k(n-k)} c_i(\Tkn)t^i \ = \ \prod_{i,j}(1 + (\mu_j -
\lambda_i)t ),
\]
where $i$ and $j$ in the product range over possible indices
\cite[{\it ibid}]{BottTu}. In other words, for each i we have
$c_i(\Tkn) = \sigma_i(\{\mu_j-\lambda_i\}_{i,j})$. This shows that
each $c_i(\Tkn)$ is symmetric in the $\lambda_i$'s and the $\mu_j$'s,
and can therefore be written in terms of the $x_i$'s and $y_j$'s.

In particular, the Euler class $e(\Gkn)$ can be lifted to a polynomial
\[
P \in K[x_1, \ldots, x_k, y_1, \ldots, y_{n-k}]
\]
or, using the relations between the $x_i$'s and $y_j$'s, to a
polynomial
\[
P' \in K[x_1, \ldots, x_k].
\]
$P'$ is referred to as the ``Euler polynomial.''

Bertram \cite{Ber94} has proven the following:
\begin{prop}
   \label{prliftHess} For each $(k,n)$, the Euler polynomial $P'$ is a
lifting of $(-1)^{n \choose 2}H_q \in QH^*(\Gkn)$.
\end{prop}

We can now prove theorem \ref{thhesseuler}.

\begin{prf} {\bf (of theorem \ref{thhesseuler})}
Proposition \ref{prHessomega} shows that $e_q:=e_q(\Gkn) = \kappa H_q$
for some $\kappa \in K$. Let $\pi \colon QH^*(\Gkn) \rightarrow
H^*(\Gkn)$ denote the module homomorphism sending $q \mapsto 0$. By
definition, $P'$ is a lifting of $e:=e(\Gkn)$, so (by
\ref{prliftHess}) we have
\[
\begin{diagram}
   \node{K[x_1, \ldots, x_k]} \arrow{s} \arrow{se} \\ \node{H^*(\Gkn)}
   \node{QH^*(\Gkn)} \arrow{w,t}{\pi}
\end{diagram} \quad
\begin{diagram}
    \node{P'} \arrow{s} \arrow{se} \\ \node{e(\Gkn)} \node{(-1)^{n
    \choose 2}H_q } \arrow{w,t}{\pi}
\end{diagram}
\]
where the unlabeled arrows are the canonical projection maps. Now
$\pi(e_q) = e$, by definition of the quantum multiplication $\ast$, so
proposition \ref{prHessomega} shows that
\[
(-1)^{n \choose 2}e \, = \, \pi \left( (-1)^{n \choose 2}e_q \right)
                      \, = \, \pi \left( (-1)^{n \choose 2}\kappa H_q
                      \right) \, = \, \kappa e,
\]
and thus $\kappa = (-1)^{n \choose 2}$.
\end{prf}

\section{Quantum Cohomology of Hyperplanes}
   \label{secqcohhyp}

For the sake of contrast with the Grassmannians, this section provides
another class of examples of a quantum cohomology ring, and determines
which of these are semisimple. In \cite{TiaXu96}, Tian and Xu discuss
a more general class of examples along these lines from the point of
view of semisimplicity in the sense of Dubrovin (as defined in the
introduction).

Let $X \subset \Cm P^{n+r}$ be a smooth complete intersection of
degree $(d_1, \ldots, d_r)$ and dimension $n \geq 2$ satisfying $n
\geq \sum(d_i-1) - 1$. Let $\Gamma$ denote the hyperplane class
generating $H^2(X,\Z)$. By the ``primitive cohomology $H^n(X)_0$ of
$X$'' we mean $H^n(X)$ if $n$ is odd, and the subspace of $H^n(X)$
orthogonal to $\Gamma^{n/2}$ if $n$ is even. Beauville shows in
\cite{Bea95} (although he unnecessarily presumes $q=1$), that
$QH^*(X)$ is the algebra over $K[q,q^{-1}]$ generated by $\Gamma$ and
$H^n(X)_0$, subject to the relations
\[
\Gamma^{n+1}\, = \, d_1^{d_1} \cdots d_r^{d_r}\Gamma^{d-1}q
\]
and, for all $a,b \in H^n(X)_0$,
\[
\Gamma a = 0 \ \ \ \mathrm{and}\ \ \ ab = \langle a,b \rangle
\frac{1}{d} \left(\Gamma^n - d_1^{d_1} \cdots d_r^{d_r}\Gamma^{d-2}q
\right).
\]
Here, $\langle \cdot, \cdot \rangle$ denotes the classical
intersection form $a \otimes b \mapsto f(a \cup b)$, where
$f:=(\Gamma^n)^*$.

\begin{prop}
Let $X$ denote a hyperplane of degree $d$. For any nonzero $r \in K$,
if $d > 2$ then $(\theta_r)_*[QH^*(X)]$ is not semisimple. If $d=2$
then $(\theta_r)_*[QH^*(X)]$ is semisimple.
\end{prop}
\begin{prf}
Denote $(H^n(X)_0:K)$ by $R$, and choose a basis $e_1, \ldots, e_R$
for $H^n(X)_0$. Together with the elements $1, \Gamma, \Gamma^2,
\ldots, \Gamma^n$, this provides a full vector-space basis for
$QH^*(X)$. Thus the characteristic element of $(QH^*(X),f)$ is
\[
\begin{array}{rcl}
   \omega & = & \sum_{i=0}^n \Gamma^i \Gamma^{n-i} + \sum_{i=1}^R
      e_ie_i^{\#} \\ \\ & = & (n+1)\Gamma^n +
      \frac{R}{d}\left(\Gamma^n - d^d \Gamma^{d-2}q \right).
\end{array}
\]
Notice that if $d > 2$ then $\omega$ is divisible by $\Gamma$, so
$\omega e_1 = 0$ (for example). Since $e_1 \notin \Ker \thetatil$ (it
is a basis element for any choice of coefficients!), proposition
\ref{prinducedomega} shows that $(\theta_r)_*[QH^*(X)]$ is not
semisimple for any choice of $r \in K$.

If $d=2$ then we have
\[
\omega \, = \, (n+1+\frac{R}{2})\Gamma^n - 2Rq,
\]
and thus $\Gamma\omega = 4(n+1)q\Gamma$. Order the basis for $QH^*(X)$
as follows: $1, \Gamma, \Gamma^2, \ldots, \Gamma^n,$ $e_1, \ldots,
e_R$. Then the matrix $[\bbar(\omega)]$ corresponding to $\omega$
under the regular representation $\bbar$ is
\[
\left( \begin{array}{ccccc|ccc} -2Rq & & & & n+1+\frac{R}{2} & \\ &
  4(n+1)q & & $\LARGE{0}$ & \\ & & 4(n+1)q & & & & \\ & $\LARGE{0}$ &
  & \ddots & & \\ & & & & 4(n+1)q & \\ \hline & & & & & -2nq & & 0 \\
  & & & & & & \ddots & \\ & & & & & 0 & & -2nq \end{array} \right) .
\]
Since the determinant of this matrix is a unit in $K[q,q^{-1}]$, we
see that $\omega$ is a unit in $QH^*(X)$; by \ref{prinducedomega} this
shows that $(\theta_r)_*[QH^*(X)]$ is semisimple for any choice of
$r$.
\end{prf}

\section*{Acknowledgements}

The author thanks Jack Morava for his guidance of the research leading
to this article. Aaron Bertram, Steve Sawin and Geng Xu were also
helpful in correspondence and conversation, and Peter Landweber
graciously assisted with proof-reading.


\noindent
Rutgers University

\noindent
labrams@math.rutgers.edu

\end{document}